\pdfoutput=1
\documentclass{amsart}
\usepackage[latin1]{inputenc}
\usepackage[english]{babel}
\usepackage{csquotes}

\usepackage{color} 
\usepackage[colorlinks=true,linktoc=all,linkcolor=blue,citecolor=red,filecolor=pink,urlcolor=blue]{hyperref}
\usepackage[hyperref=true,style=alphabetic,citestyle=alphabetic,backend=bibtex,maxnames=100,doi=false,isbn=false,url=false,giveninits=true]{biblatex}
\bibliography{biblio}

\usepackage{graphicx}
\usepackage{pgf,tikz}
\usetikzlibrary{arrows}
\usetikzlibrary{decorations.markings}

\usepackage{amssymb,amsmath,latexsym}
\theoremstyle{plain}                    
\newtheorem{theorem}{Theorem}[section]
\newtheorem{lemma}[theorem]{Lemma}
\newtheorem{proposition}[theorem]{Proposition}
\newtheorem{corollary}[theorem]{Corollary}
\theoremstyle{definition}
\newtheorem{definition}[theorem]{Definition}
\newtheorem{example}[theorem]{Example}

\theoremstyle{remark}
\newtheorem{remark}[theorem]{Remark}

\numberwithin{equation}{section}

\newcommand{\real}{\mathbb{R}}
\newcommand{\cmp}{\mathbb{C}}
\newcommand{\intz}{\mathbb{Z}}
\newcommand{\nat}{\mathbb{N}}
\newcommand{\rat}{\mathbb{Q}}
\newcommand{\hyp}{\mathbb{H}}

\newcommand{\cp}{\mathbb{C}\mathbb{P}^1}
\newcommand{\pslc}{\mathrm{PSL}_2\cmp}
\newcommand{\pslr}{\mathrm{PSL}_2\real}

\newcommand{\psls}{\text{M\"ob}(\sigma)}
\newcommand{\pslspar}[1]{\text{M\"ob}(#1)}
\newcommand{\pslss}{\text{\emph{M\"ob}}(\sigma)}
\newcommand{\pslsspar}[1]{\text{\emph{M\"ob}}(#1)}
\newcommand{\aut}{\text{Aut}}
\newcommand{\aaut}{\text{\emph{Aut}}}

\begin{document}

\title[$\cp$-structures with maximal number of M\"obius transformations]{Complex projective structures with maximal number of M\"obius transformations}

\author[G. Faraco]{Gianluca Faraco}
\address{Dipartimento di Matematica - Universit\`a di Parma, Parco Area delle Science 53/A, 43132, Parma, Italy}
\curraddr{School of Mathematics - Tata Institute of Fundamental Research, Homi Bhabha Road, Navy Nagar, Colaba, Mumbai 400005, India}
\email{gianluca.faraco.math@gmail.com}

 \author[L. Ruffoni]{Lorenzo Ruffoni}
\address{Dipartimento di Matematica - Universit\`a di Bologna, Piazza di Porta San Donato 5, 40126, Bologna, Italy}
\curraddr{Department of Mathematics - Florida State University, 1017 Academic Way, Tallahassee, FL 32306-4510, USA}
\email{lorenzo.ruffoni2@gmail.com}

\keywords{Complex projective structures, Fuchsian groups, automorphisms, Galois Bely\u{\i} curves, quasiplatonic surfaces}
\subjclass[2010]{57M50,14H37}%
\date{\today}

\begin{abstract}
  We consider complex projective structures on Riemann surfaces and their groups 
of projective automorphisms. We show that the structures achieving the maximal 
possible number of projective automorphisms allowed by their genus are 
precisely the Fuchsian uniformizations of Hurwitz surfaces by hyperbolic 
metrics. More generally we show that Galois Bely\u{\i} curves are precisely those 
Riemann surfaces for which the Fuchsian uniformization is the unique complex 
projective structure invariant under the full group of biholomorphisms.
\end{abstract}

\maketitle                
\tableofcontents

\section{Introduction}
Upgrading a smooth surface $S$ to a Riemann surface $X$ by 
introducing a complex structure turns it into a quite rigid object from several 
points of view. The prototypical phenomenon was discovered by Hurwitz in 1893, 
who proved in \cite{HU} that, if the genus $g$ of the surface is at least 2, 
then the group of holomorphic automorphisms of $X$ is finite, the 
cardinality being bounded just in terms of the genus as $\#\aut(X)\leq 
84(g-1)$. Riemann surfaces which attain this bound have 
been named Hurwitz surfaces and the finite groups arising as their groups of 
automorphisms have been 
named Hurwitz groups. A lot of research has then been done to study geometric and algebraic 
properties of these objects: see \cite{CO2}, \cite{CO}, \cite{LE} and the references 
therein.\par
  It is natural to ask what happens if we further rigidify the geometry, for 
instance fixing a special coordinate covering for a Riemann surface 
(in the sense of \cite{GU}), i.e. an open cover for which transition functions 
are not just local biholomorphisms, but belong to some more
restricted group of transformations of the Riemann sphere $\cp$. 
Schlage-Puchta and Weitze-Schmith\"usen in \cite{SPWS}
consider the case of translation surfaces, for which transition functions are translations. They 
show that a result similar to Hurwitz's one 
holds, with an explicit bound being $4(g-1)$; they call Hurwitz translation surfaces the structures 
which achieve this bound, and characterise them as normal origamis, i.e. 
square-tiled surfaces which arise as Galois coverings of the standard torus 
$\cmp / \intz [i]$.\par
  Here we consider the same questions for complex  projective 
structures, i.e. geometric structures defined by an atlas of charts for which 
transitions are given by restrictions of global 
M\"obius transformations of the Riemann sphere. Such a 
structure in particular induces a complex structure on the underlying surface, 
and indeed the interest in this kind of structures comes from classical 
uniformization theory (see \cite{GU}) as well as from the study of linear ODEs 
(see \cite{GKM}). We show that the group of projective automorphisms of a 
complex projective structure can be as large as the group of 
holomorphic automorphisms of the underlying Riemann surface, and that the
complex projective structures which attain this bound are precisely the ones arising from Fuchsian
uniformizations of Hurwitz surfaces by hyperbolic metrics.\par
  The classical approach to this kind of questions is to consider  the 
quotient of a structure by its full group of automorphisms; however this does 
not work for projective structures, since the 
quotient might not exist as a projective structure. Our approach is 
to take a relative point of view: given a Riemann surface, we ask which are the complex projective 
structures on it for which the full
group of holomorphic automorphisms acts by projective automorphisms; we call 
such structures relatively Hurwitz projective structures. For genus $g\geq 2$ Fuchsian 
uniformizations are always among them, but generically every structure has 
this property, just because a generic Riemann surface of genus $g\geq 3$ has no non-trivial 
automorphism. We identify the condition under which the Fuchsian uniformization 
is the only structure with such a property; our main result (see Theorem 
\ref{mainthm} below) can be stated as follows.

\begin{theorem}
 Let $X$ be a Riemann surface of genus $g\geq 2$. Then the following are 
equivalent:
 \begin{itemize}
  \item[(i)] $X$  is a Galois Bely\u{\i} curve; 
  \item[(ii)] the Fuchsian uniformization of $X$ is the unique relatively Hurwitz projective structure in $\mathcal{P}(X)$.
 \end{itemize}
\end{theorem}

  An analogous statement holds in genus $g=1$ replacing ``Fuchsian'' by ``Euclidean''. The first condition in the above statement means that 
$X/\aut(X)$ has genus zero and that $\pi:X\to X/\aut(X)$ branches exactly over 3 points, so that $\pi$ is both a Bely\u{\i} function on $X$ and a Galois cover; in particular these Riemann surfaces are examples of Bely\u{\i} curves, i.e. algebraic curves which can be defined over the field of algebraic numbers $\overline{\rat}$. These curves are also known as quasiplatonic surfaces, and have been studied by many authors (see for instance \cite{CV}, 
\cite{JW}, \cite{POP}, \cite{SPW} and \cite{WO});  they are known to be exactly those Riemann surfaces which have strictly more automorphisms than any proper deformations, or equivalently that can be uniformized by normal subgroups of hyperbolic triangle groups. A straightforward consequence of the theorem is that a complex projective structure has strictly more automorphisms than any proper deformation if and only if it is the Fuchsian uniformization of a Galois Bely\u{\i} curve; in particular we recover the aforementioned  characterization  of Hurwitz projective structures as Fuchsian uniformizations of Hurwitz Riemann surfaces (see Corollary \ref{maincor} below).\par 
  As in the case of translation surfaces, 
we have therefore a very neat geometric description of the structures which 
attain the bound. Moreover in both cases these structures turn out to be 
integral points of their moduli spaces: normal origamis are among square-tiled 
surfaces, which are integral points in period coordinates (see \cite{Z}), and 
Fuchsian uniformizations are among Fuchsian projective structures, which are 
integral points in Thurston's coordinates  (see \cite{DU}).\par

  The paper  is organised as follows: Section 2 contains the basics about complex projective 
structures and their automorphisms, and the definition of (relatively) Hurwitz projective structure; 
in Section 3 we describe the picture for genus 1 and for translation surfaces. In Section 4 we 
gather the technical lemmas about the Schwarzian parametrization of projective structures which are 
used in Section 5 to prove the main result; in particular a criterion is obtained for a 
biholomorphism $f\in \aut(X)$ to act projectively on a given projective structure in terms of an 
affine action of $\aut(X)$ on the space of holomorphic quadratic differentials $H^0(X,K^2)$.\par
\vspace{.5cm}
\textbf{Acknowledgements}:
We would like to thank Stefano Francaviglia for his encouragement in developing this project, Pascal Hubert for pointing us to useful references, and the referees for carefully reading this paper. The second author acknowledges support by INDAM-GNSAGA.

\section{Complex projective structures and their automorphisms}
  In  this section we gather the needed preliminaries on 
complex projective structures that we will use in the sequel. 

\subsection{Basic definitions} Let $S$ be a closed, connected and 
orientable surface, let $\cp=\cmp \cup \{\infty\}$ be the Riemann sphere and let $\pslc$ be its 
group of holomorphic automorphisms acting by M\"obius transformations
 \[ \pslc \times \cp \to \cp, \left(\begin{array}{cc} 
 a & b\\ c & d \\ \end{array} \right),z \mapsto \dfrac{az+b}{cz+d} \]

\begin{definition} A complex projective structure  $\sigma$ on $S$ is 
a maximal atlas of charts taking values in $\cp$ and such that transition 
functions are 
restrictions of elements in $\pslc$. 
\end{definition}

  In the following we will also refer to such structures simply as 
projective structures. Given a projective structure $\sigma$, performing 
analytic continuation of local 
charts of $\sigma$ along paths in $S$ gives rise to an immersion $dev:\widetilde{S}\to \cp$, 
usually 
called a developing map for $\sigma$, well-defined up to post-composition by a M\"obius 
transformation.
We refer to \cite{DU} for a detailed survey about complex projective structures.

\begin{remark}\label{rmk_unif}
Since M\"obius transformations are in particular holomorphic maps, a projective structure 
on $S$ always determines an underlying Riemann surface structure on $S$. Conversely, by the 
classical uniformization theory, any Riemann surface $X$ is of the form $U/\Gamma$ 
where $U$ is an open subset of $\cp$ and $\Gamma$ is a discrete subgroup of $\pslc$ acting freely 
and properly discontinuously on $U$; this endows $X$ with a complex projective structure, namely 
the one coming from the identification $X\cong U/\Gamma$. When $g=0$ there is only one possible 
choice $U=\cp$ and $\Gamma=1$. If $g=1$ then $U$ can be chosen to be the complex plane $\cmp$ and 
$\Gamma$ to be a discrete free abelian group of rank 2 acting by translations. When $g\geq 2$ it is 
possible to choose $U$ to be the upper-half plane $\mathcal{H}=\{z \in \cmp \ | \ Im(z)>0\}$ and 
$\Gamma$ to be a Fuchsian group, i.e. a discrete subgroup of $\pslr$. Notice in particular that for 
a projective structure of this type any developing map is a diffeomorphism $dev:\widetilde{S}\to 
U$, and that it endows $X$ with a canonical Riemannian metric of constant sectional curvature 
$k=sign(\chi(X))\in \{1,0,-1\}$.
\end{remark}

\begin{definition}\label{def_unif}
 Let $X$ be a Riemann surface of genus $g$. If $\Gamma \subset \pslc$ is discrete, $U\subset \cp$ 
is an open connected domain on which $\Gamma$ acts freely and properly 
discontinuously with quotient $U/\Gamma$ diffeomorphic to $S$, then we say 
that the projective structure $\sigma=U/\Gamma$ is a uniformization of $X$. In particular if $g=1$ 
then we call the Euclidean uniformization of $X$ the complex projective structure coming 
from 
the flat Riemannian metric as in Remark \ref{rmk_unif}, and if $g\geq 2$ then we call the 
Fuchsian uniformization of $X$ the complex projective structure coming from the hyperbolic 
Riemannian metric as in Remark \ref{rmk_unif}.
\end{definition}

\begin{remark}
Among all possible ways of realising a Riemann surface $X$ as a quotient $U/\Gamma$ as in Remark 
\ref{rmk_unif}, Euclidean and Fuchsian uniformizations are special in that $U$ is simply 
connected and $\Gamma$ is isomorphic to the fundamental group $\pi_1(X)$; other uniformizations 
arise for instance from Schottky uniformization, where $U$ is the complement of a Cantor set and 
$\Gamma$ is a (non-abelian) free group. On the other hand, not every complex projective structure 
is of the form $U/\Gamma$: for instance in \cite{MAS2} Maskit has produced many examples of 
projective structures with surjective and non injective developing maps.
\end{remark}

  Now we turn to the study of maps between projective structures.
\begin{definition}  
Let $\sigma_1,\sigma_2$ be projective structures on $S$ and let $f:\sigma_1\to \sigma_2$ be a 
diffeomorphism. We say that $f$ is  projective if its restrictions to 
local  projective charts are given by elements in $\pslc$. We say that 
$\sigma_1$ and $\sigma_2$ are 
isomorphic if there exists a projective diffeomorphism between them.
\end{definition}
\begin{definition}
 Let $\sigma$ be a projective structure. We define the group of projective automorphisms 
of $\sigma$ to be $\psls=\{ f \in \textsf{Diff}(S) \ | \ f \textrm{ is projective for } \sigma \}$.
\end{definition}

  Up to isomorphism there is a unique complex projective structure $\sigma$ on a surface of 
genus 0, given by its realisation as $\cp$; of course  $\psls=\pslc$. The following straightforward 
observations will be useful in the following.
\begin{lemma}\label{lem_projviadev}
 Let $\sigma$ be a projective structure on $S$ and $dev:\widetilde{S}\to \cp$ be a developing map 
for $\sigma$. Then $f \in \emph{\textsf{Diff}}(S)$ is projective if and only if there 
exists $g \in 
\pslc$ such that $dev\circ \widetilde{f}=g\circ dev$, for some lift $ \widetilde{f}$ of $f$ to  
$\widetilde{S}$.
\end{lemma}

  Let us denote by $\aut(X)$ the group of holomorphic automorphisms of a Riemann surface 
$X$. The following is immediate.

\begin{lemma}\label{lem_projisholo}
Let $\sigma$ be a projective structure on $S$ and let $f\in \pslss$. If $X$ is the underlying 
Riemann surface, then $f \in \aaut(X)$.\\
\end{lemma}

\subsection{Hurwitz projective structures}
From now on, unless otherwise specified, we assume that the surface $S$ has genus $g\geq 2$. 
Let $\sigma$ be a projective structure on $S$ and let $X$ be the underlying 
Riemann surface. By the classical Hurwitz Theorem the group $\aut(X)$ of 
holomorphic automorphism of $X$ is a finite group of cardinality $\# 
\text{Aut}(X)\leq 84(g-1)$ (see \cite[V.1.3]{FK}). A Hurwitz surface is 
classically defined to be a Riemann surface of genus $g$ whose automorphism group attains this 
bound. By the above  Lemma  \ref{lem_projisholo}  the same bound 
holds for the group $\psls$ of $\sigma$, and one 
can ask if it can be sharpened or not.
  A straightforward approach is to pick a Hurwitz Riemann 
surface $X$ and look for projective structures for which 
holomorphic automorphisms are also projective. This might be 
 hard to check in general, but the collection of uniformizations of $X$ (see 

Definition

 \ref{def_unif}) provides a playground where we can perform concrete 
computations.
\begin{lemma}\label{lem_covering}
 Let $X$ be a Riemann surface and $\sigma=U/\Gamma$ a uniformization of $X$. Suppose that every 
$f\in \aaut(X)$ lifts to some $\widetilde{f}\in \aaut(U)$ and that $\widetilde{f}$ 
is the restriction of a M\"obius transformation. Then $\pslss=\aaut(X)$.
\end{lemma}
\begin{proof}
By classical covering space theory, lifts of automorphisms of $X$ are exactly 
given by automorphisms 
of $U$ which normalise $\Gamma$ (both in the holomorphic and in the projective setting); in other 
words the normalizer $N_{\aut(U)}(\Gamma)$ of $\Gamma$ in $\aut(U)$ acts by biholomorphisms on $X$, 
the normalizer $N_{\pslspar{U}}(\Gamma)$ of $\Gamma$ in $\pslspar{U}$ acts  by projective 
diffeomorphisms on $\sigma$ and both actions factor exactly through the standard action of $\Gamma$ 
by deck transformations of the covering $U\to U/\Gamma$; so we get injective maps 
$N_{\aut(U)}(\Gamma)/\Gamma \to \aut(X)$ and $N_{\pslspar{U}}(\Gamma)/\Gamma\to \psls$. 
The hypothesis imply that $N_{\aut(U)}(\Gamma)/\Gamma \to \aut(X)$ is surjective and that 
$N_{\aut(U)}(\Gamma)/\Gamma=N_{\pslspar{U}}(\Gamma)/\Gamma$. In particular we obtain an injective 
map $\aut(X)\to \psls$. Since by  Lemma  \ref{lem_projisholo} we always have $\psls \subseteq \aut(X)$ and 
both groups are finite, we get the result.
\end{proof}
  We can use the previous lemma to obtain that in general the Hurwitz 
bound can not be improved.
\begin{lemma}\label{lem_fuchshurwishurwproj}
 Let $X$ be a Hurwitz surface and $\sigma$ its Fuchsian uniformization. Then 
$\pslss=\aaut(X)$; in particular $\#\pslss = 84(g-1)$.
\end{lemma}
\begin{proof}
By definition of Fuchsian uniformization, $\sigma$ is of the form $\mathcal{H}/\Gamma$, where 
$\mathcal{H}$ is the upper-half plane and $\Gamma \subset \pslr$ is a Fuchsian group.
In particular $\mathcal{H}$ is simply connected and it is well-known that its biholomorphisms are 
exactly given by M\"obius transformations with real coefficients, i.e. 
$\aut(\mathcal{H})=\pslspar{\mathcal{H}}=\pslr$; therefore we can just apply  Lemma 
\ref{lem_covering}.
\end{proof}
  We give therefore the following definition, by analogy with the 
classical case.
\begin{definition}
 Let $\sigma$ be a projective structure on $S$. We say that $\sigma$ is a 
Hurwitz projective structure if $\# \psls = 84(g-1)$.
\end{definition}

  Lemma \ref{lem_fuchshurwishurwproj} can then be restated by saying 
that Fuchsian uniformizations of Hurwitz surfaces are Hurwitz projective 
structures. Conversely one can ask if projective 
structures with a maximal number of projective automorphisms are necessarily 
given by hyperbolic metrics on Hurwitz surfaces. The answer turns out to be 
positive, and we will recover this result as consequence of a more general 
statement below (see  Corollary  \ref{maincor}).

\begin{remark}
The properties of Fuchsian uniformizations which we have used in the proof of  Lemma 
\ref{lem_fuchshurwishurwproj}  to apply  Lemma  \ref{lem_covering} are the following: the domain $U$ is 
simply connected and $\aut(U)=\pslspar{U}$. By \cite{MIN} such properties are 
known to characterise Fuchsian uniformizations. The same kind of arguments might fail for other 
uniformizations. For instance if $U$ is an open quadrant then it is simply connected, but not every 
biholomorphism is the restriction of a M\"obius transformation. On the other hand if $U$ is a 
Schottky domain then it is not simply connected, but every biholomorphism is indeed projective (by 
\cite[IV.2D,IV.19F]{AS}).
\end{remark}

\begin{remark}\label{hownotto1}
The study of geometric automorphisms of surfaces is usually done by taking the quotient by the action of the automorphism group, and studying the geometry and combinatorics of the resulting quotient map.  In general this action may have fixed points, so that the quotient is endowed with an orbifold geometric structure.
In the case of Riemann surfaces it turns out that the orbifold charts can actually be uniformized in such a way that the quotient carries a honest complex structure with respect to which the projection is holomorphic (see \cite[Theorem III.3.4]{MIR}). Translation surfaces are more rigid, but the key observation  (see \cite{SPWS}) is that fixed points of a translation can appear only at cone points of the structure, so that the effect of taking the quotient is absorbed in the reduction of the total curvature, and also in this case the quotient is endowed with a honest translation structure. In both cases this is quite useful in taming the behaviour of the quotient map (e.g. its degree).
On the other hand, for complex projective structures it is in general not possible to reuniformize the orbifold charts in such a way that the quotient supports a honest projective structure with respect to which the quotient map is projective. Consider for instance a closed  surface of genus $g\geq 2$  endowed with a hyperelliptic hyperbolic structure (seen as a projective structure here); the structure induced on the quotient is a hyperbolic cone structure on a genus 0 surface, which is of course not the standard projective structure on $\cp$.
We will bypass this issue by looking at the action of the automorphism group on the space of holomorphic quadratic differentials, which is classically known to provide a parametrization of the space of projective structures (see \ref{def_schwarzpar} below). The previous comments translate into this setting by saying that a holomorphic quadratic differential descends to a meromorphic, but in general not holomorphic, differential on the quotient (see Remark \ref{hownotto2} below).
\end{remark}

\subsection{Relatively Hurwitz projective structures}
Our approach to the problem will be to consider it in a relative way: any projective structure has 
an underlying complex one, and by  Lemma  \ref{lem_projisholo} any projective automorphism is in particular 
holomorphic for it; a natural question is to ask when every holomorphic automorphism is indeed 
projective, which motivates the following definition.
\begin{definition}
 Let $\sigma$ be a projective structure on $S$ and $X$ the underlying Riemann surface. We say that 
$\sigma$ is a relatively Hurwitz projective structure if $ \psls = \aut(X)$.
\end{definition}
\begin{example}\label{ex_trivialrelHurproj}
 A generic Riemann surface $X$ of genus $g\geq 3$ has no non-trivial automorphisms, therefore any 
projective structure on $X$ is a relatively Hurwitz projective structure. 
\end{example}

\begin{example}
 On the other end of the 
spectrum, let $X$ be a Hurwitz surface. By  Lemma  \ref{lem_fuchshurwishurwproj} the Fuchsian 
uniformization $\sigma$ of  $X$ is a relatively Hurwitz projective structure: any 
holomorphic automorphism is indeed an isometry for the uniformizing hyperbolic metric, and those 
are all projective.  
\end{example}

  More generally, we can apply  Lemma  \ref{lem_covering} as in  Lemma  \ref{lem_fuchshurwishurwproj} to 
obtain the following.
\begin{lemma}\label{lem_fuchsunifisrelhurwproj}
 Let $X$ be any Riemann surface and $\sigma$ its Fuchsian uniformization. Then 
$\pslss=\aaut(X)$; in particular $\sigma$ is a relatively Hurwitz projective structure.
\end{lemma}

\begin{example}
 Non-trivial examples of relatively and non-relatively Hurwitz projective structures will be given 
below in  Example  \ref{ex_relHurproj} on hyperelliptic Riemann surfaces of genus 
$g\geq 3$ whose automorphism group is generated by the hyperelliptic involution.
\end{example}

  The following  question turns out to hide rich geometric phenomena: when is it true that 
the Fuchsian uniformization of $X$ is the unique relatively Hurwitz projective structure? As a 
warm-up, we will first consider the analogous problem in other settings, namely translation 
surfaces and elliptic curves.

\section{An affine detour}
  In this section we gather some results available in the literature in order to 
describe the geometric features enjoyed by geometric structures with maximal 
group of automorphisms in the setting of affine geometry.

\subsection{Complex affine structures on tori}\label{section_tori}
  In this section (and only in this section) we consider the situation in which $S$ has 
genus $g=1$. In this case complex projective structures are actually complex affine structures (see 
\cite{GU}), i.e. the local charts take values in a fixed affine patch $\cmp\subset \cp$ and the 
changes of coordinates are restrictions of complex affine transformations $g\in 
\textsf{Aff}_1(\cmp)=\{z\mapsto az+b \ | \ a \in \cmp^*, b \in \cmp\}$. These geometric structures 
are conveniently described as affine deformations of the Euclidean uniformizations (see \cite{GU} or 
\cite{LM} for more details). More precisely let $\tau \in 
\mathcal{H}^+/\mathrm{SL}_2\intz=\{Im(z)>0\}/\mathrm{SL}_2\intz$, 
$\Lambda_\tau=span_\intz\{1,\tau\}$ and $X_\tau=\cmp/\Lambda_\tau$. The complex affine structures on 
$X_\tau$ are parametrized by $c \in \cmp$: for $c=0$ we have the Euclidean uniformization $\sigma_E$ 
of $X_\tau$ (see  Remark  \ref{rmk_unif}); for $c \in \cmp^*$ we get a non-Riemannian affine structure 
$\sigma_c$ on $X_\tau$ defined by a developing map $dev_c:\cmp\to \cmp, z \mapsto e^{cz}$. The group 
of automorphisms $\aut(X_\tau)$ of the torus $X_\tau$ is the semidirect product $\aut(X_\tau)=X_\tau 
\rtimes \aut^0(X_\tau)$, where $X_\tau$ acts on itself by translations $T_p(z)=z+p$ for $ p\in 
X_\tau$ and $\aut^0(X_\tau)$ is a finite group which generically consists only of the hyperelliptic 
involution $J(z)=-z$, and also contain a complex  multiplication $R_\tau(z)=\tau z$,when $\tau = 
e^{i\frac{2\pi}{3}}$ or $\tau=e^{i\frac{\pi}{2}}$. A straightforward computation shows that the 
translation part always acts by complex affine transformations: this is clear for the Euclidean 
uniformization, and for the other structures we just observe that a translation acts 
 as  $dev_c(T_p(z))=dev_c(z+p)=e^{c(z+p)}=e^{cp}dev_c(z)$, i.e. as a complex dilation, and use  Lemma 
\ref{lem_projviadev}.

\begin{definition}
Let $X$ be a Riemann surface of genus 1. We denote by $\mathcal{A}(X)$ the space of complex affine 
structures on $X$. If $\sigma \in \mathcal{A}(X)$, then we denote by $\textsf{Aff}_1(\sigma)$  the 
group of affine automorphisms, i.e. the ones that are given by elements in $\mathsf{Aff}_1(\cmp)$ in 
affine charts for $\sigma$, and by $\textsf{Aff}_1^0(\sigma)$ the reduced group obtained by 
quotienting out the translations, so that $\textsf{Aff}_1(\sigma)=X_\tau \rtimes 
\textsf{Aff}_1^0(\sigma)$.
\end{definition}

  Notice that in general we only have $\textsf{Aff}_1^0(\sigma)\subseteq\aut^0(X)$. In 
particular $\textsf{Aff}_1^0(\sigma)$ has at most order 6, which is the bound for the cardinality of 
the reduced automorphism group $\aut ^0(X)$ of a Riemann surface of genus 1.

\begin{example}\label{ex_euclidean}
A direct computation shows that actually each automorphism of $X$ is a complex affine transformation 
(indeed a Euclidean isometry) with respect to its Euclidean uniformization $\sigma_E$, i.e. 
$\aut(X)=\mathsf{Aff}_1(\sigma_E)$. In particular  $\textsf{Aff}_1^0(\sigma_E)=\aut^0(X)$; therefore 
the Euclidean uniformization of the complex torus $X_{e^{i\frac{2\pi}{3}}}$ achieves the maximum 
possible cardinality.
\end{example}

  We are naturally led to give the following definitions.

\begin{definition}
Let $X$ be a Riemann surface of genus 1 and $\sigma \in \mathcal{A}(X)$; we say $\sigma$ is a 
Hurwitz affine structure if $\#\textsf{Aff}_1^0(\sigma)=6$, and that it is a relatively 
Hurwitz affine structure if $\textsf{Aff}_1^0(\sigma)=\aut^0(X)$.
\end{definition}

  We then have the following.

\begin{proposition}\label{prop_affine}
The unique Hurwitz affine structure is the Euclidean uniformization of $X_{e^{i\frac{2\pi}{3}}}$. 
For a Riemann surface $X$ of genus 1 the unique relatively Hurwitz affine structure in $ 
\mathcal{A}(X)$ is its Euclidean uniformization.
\end{proposition}

\begin{proof} 
The first statement follows  from the above discussion and the fact that $X_{e^{i\frac{2\pi}{3}}}$ 
is the only torus with $\aut^0(X_{e^{i\frac{2\pi}{3}}})=6$. For the second statement fix $X$, let 
$J\in \aut(X)$ be the hyperelliptic involution and assume $\sigma$ is not the Euclidean 
uniformization. Then we just compute that $J$ acts as 
$dev_c(J(z))=dev_c(-z)=e^{-cz}=\frac{1}{dev_c(z)}$, i.e. as an inversion, which is not affine. 
\end{proof}

  Looking at affine structures from a projective point of view, we see that the Euclidean 
uniformization of $X_{e^{i\frac{2\pi}{3}}}$ is a Hurwitz projective structure (by  Example 
\ref{ex_euclidean}). The same ideas as above also prove the following.

\begin{proposition}
  Let $X_\tau$ be a Riemann surface of genus 1. Then the following are equivalent
\begin{itemize}
\item[(i)]  $\tau = e^{i\frac{2\pi}{3}}$ or $\tau=e^{i\frac{\pi}{2}}$
\item[(ii)] the Euclidean uniformization of $X_\tau$ is the unique relatively Hurwitz 
projective structure in $\mathcal{A}(X_\tau)$.
\end{itemize}
\end{proposition}

\begin{proof}
From the previous proof of  Proposition  \ref{prop_affine} we know that the hyperelliptic involution does not act 
affinely, but at least it always acts projectively. But if we consider surfaces which admit complex 
multiplication, we see that $R_\tau \in \aut(X_\tau)$ acts on a structure which is not the Euclidean 
one as $dev_c(R_\tau(z))=dev_c(\tau z)=e^{\tau cz}=dev_c(z)^\tau$, i.e. not projectively.
\end{proof}

  We see therefore that requiring the Euclidean uniformization to be the unique relatively 
Hurwitz projective structure picks out special points in the moduli space of genus 1 Riemann 
surfaces, namely those having extra automorphisms, or equivalently, those corresponding to regular 
lattices. In the following sections we are going to extend this result to higher genus surfaces, 
replacing Euclidean uniformization by Fuchsian uniformization.

\subsection{Translation surfaces}\label{section_trans}
  We now turn to consider translations surfaces. These can be defined 
either as couples $(X,\omega)$ where $X$ is a Riemann surface and $\omega \in 
H^0(X,K),\omega \neq 0$, or as complex projective structures with change of 
coordinates given by translations and with branch points (corresponding to the 
zeroes of $\omega$); see \cite{Z} for more details.\par
  The problem of counting translation automorphisms of a translation 
surface, and of describing the geometric features of structures maximizing this 
number, has been considered in \cite{SPWS}, where the following definition is 
introduced.
\begin{definition}
 Let $(X,\omega)$ be a translation surface. Denote by 
$\text{Trans}(X,\omega)$ the group of automorphisms of $X$ which are given 
by translations in local charts for the translation structure.
\end{definition}
  Then the following result is obtained in \cite{SPWS}.
\begin{theorem}[Schlage-Puchta, Weitze-Schmith\"usen \cite{SPWS}]\label{spws}
 Let $(X,\omega)$ be a translation surface. Then 
$\#\emph{Trans}(X,\omega)$ $\leq  4(g-1)$. Moreover 
$\#\emph{Trans}(X,\omega)= 4(g-1)$ (i.e. $(X,\omega)$ is a 
Hurwitz translation surface) if and only if $(X,\omega)$ is a normal 
origami.
\end{theorem}
  Recall that origamis  (also known as square-tiled surfaces) are a very special type of 
translation surfaces; they can be defined as those $(X,\omega)$ arising as a covering of the 
standard torus $\cmp/\intz[i]$ branched exactly at one point, and the normal ones are those for 
which this covering is normal. The moduli space of translation surfaces carries a natural action by 
GL$_2\real$: given a translation structure $(X,\omega)$ and a matrix $A \in $GL$_2\real$ we can 
obtain a new translation structure $A.(X,\omega)$ by just multiplying by $A$ every chart from a 
translation atlas for $(X,\omega)$;  in general this will change the underlying complex structure 
unless $A\in$ GL$_1\cmp <$ GL$_2\real$ (see \cite{Z} for more details). Notice however that 
$(X,\omega)$ and $A.(X,\omega)$ have the same translation automorphisms: if a diffeomorphism 
$f\in \textsf{Diff}(S)$ is expressed as $f(z)=z+v$ in local charts for $(X,\omega)$ then it is 
given by $f(w)=w+Av$ in the corresponding local charts for $A.(X,\omega)$. Therefore the above 
statement should be considered up to this GL$_2\real$ action: Hurwitz translation surfaces are 
those whose GL$_2\real$-orbit contains a normal cover of the standard torus $\cmp/\intz[i]$ 
branched over a single point. 
For the sake of completeness,  we now discuss what would happen in 
the context of translation surfaces when taking a relative point of view analogous to the one we 
are 
adopting here. 

\begin{definition}
A translation surface $(X,\omega)$ is a relatively Hurwitz translation 
surface if $\aut(X)=\text{Trans}(X,\omega)$.
\end{definition}
  If looking for Hurwitz objects has led to the consideration of Fuchsian/Euclidean 
uniformizations in the previous setting (see  Lemma  \ref{lem_fuchshurwishurwproj} and  Example  \ref{ex_euclidean}), looking at Hurwitz translation surfaces suggests to look for the property of being a normal 
origami. One could ask, as before, if there is an interesting relation between this geometric feature and 
the relative Hurwitz condition. The main difference with respect to the previous setting is that, by 
the above bound in  Theorem  \ref{spws}, relatively Hurwitz translation structures simply do not exist on 
Riemann surfaces which have too many automorphisms, and normal origamis do not exist on Riemann surfaces 
which do not have enough.
 On the other hand, if $\# \aut(X) = 4(g-1)$, then it is 
straightforward to see that relatively Hurwitz translation surfaces on $X$ are precisely given by the normal origamis on it (up to 
GL$_2\real$-action), just by checking cardinalities and applying  Theorem \ref{spws}. 
\begin{remark}
Given a Riemann surface $X$, relatively Hurwitz translation surfaces $(X,\omega)$ are precisely 
defined by non-zero fixed points for the natural linear action of $\aut(X)$ on $H^0(X,K)$, 
because translation automorphisms are exactly those which preserve the abelian differential.
By \cite[Corollary V.2.2]{FK} the invariant differentials form a subspace of dimension equal to the 
genus of the quotient $X/\aut(X)$.
When $\# \aut(X) = 4(g-1)$, this quotient can have genus either 0 or 1; a relatively 
Hurwitz translation structure exists on $X$ precisely when it is 1 and $X\to X/\aut(X)$ is ramified 
over only one point.
In particular if such a structure $(X,\omega)$ exists over $X$, then there is precisely a 
1-dimensional family of them: this is precisely the orbit of $(X,\omega)$ under GL$_1\cmp <$ 
GL$_2\real$. Concretely, these structures are obtained by deforming it only by rotations and 
dilatations; this is the same as pulling back the differential $\lambda dz$ on the quotient torus 
via the covering $X\to X/\aut(X)$, varying the parameter $\lambda \in \cmp^*$. A necessary 
condition 
for a normal origami to exist on a surface of genus $g$ was identified by group-theoretic 
techniques 
in \cite{SPWS}, namely $g-1$ must be divisible by 2 or 3.
\end{remark}

\section{Action of biholomorphisms on projective structures}
  In this section we review the classical Schwarzian parametrization for complex projective 
structures with a fixed underlying complex structure (in the sense of  Remark  \ref{rmk_unif}), in order to 
fix terminology and notation.

\begin{definition}
Let $X$ be a Riemann surface. We denote by $\mathcal{P}(X)$ the set of projective structures whose 
underlying complex structure is $X$.
\end{definition}

  We will prove a criterion for a biholomorphism to be projective for a given projective 
structure in terms of its action on the space of holomorphic quadratic differentials.

\subsection{Schwarzian parametrization}
  The classical parametrization of $\mathcal{P}(X)$ by holomorphic quadratic differentials 
(see \cite{DU}, \cite{GU}) is achieved by means of the following differential operator.

\begin{definition}
 Let $\Omega\subset \cmp$ be an open domain and $f:\Omega \to \cmp$ be a holomorphic and locally 
injective function. The Schwarzian 
derivative of $f$ is defined to be 
\[ \mathcal{S}(f)=\left(\frac{f''}{f'}\right)'-\frac{1}{2}\left( \frac{f''}{f'}\right)^2 \]
\end{definition}

  The basic and well-known properties of this operator are the following.

\begin{lemma}\label{schwarzianprop}
Let $\Omega\subset \cmp$  be open and $f,g:\Omega \to \cmp$ be holomorphic functions such that 
$g(\Omega)\subset \Omega$. Then $\mathcal{S}(f\circ g)=( \mathcal{S}(f)\circ g )\cdot 
g'^2+\mathcal{S}(g)$. 
\end{lemma}

\begin{lemma}\label{schwarzianprop2}
Let $\Omega\subset \cmp$  be open and $f:\Omega \to \cmp$ be holomorphic function. Then 
$\mathcal{S}(f)=0$ if and only if $f$ is the restriction of a M\"obius transformation.
\end{lemma}

  A direct consequence is the following: given a projective structure $\sigma_0$ over a 
Riemann surface $X$ and a holomorphic map $f:\widetilde{X}\to \cp$ we can lift the projective 
structure from $X$ to $\widetilde X$ and we can compute the Schwarzian derivative of $f$ in local 
projective coordinates on $\widetilde{X}$. By the above formulae, we get a collection of locally 
defined functions on $\widetilde X$, which pack together to define a quadratic differential on 
$\widetilde{X}$; this is actually invariant under covering transformations, and thus descends to a 
well-defined quadratic differential on $X$, which we denote $\mathcal{S}_{\sigma_0}(f)$. In 
particular we can give the following definition.

\begin{definition}\label{def_schwarzpar}
Let $\sigma_0$ and $\sigma$ be projective structures on $X$. The Schwarzian derivative of $\sigma$ 
with respect to $\sigma_0$ is the holomorphic quadratic differential on $X$ given 
 by  $\mathcal{S}_{\sigma_0}(\sigma)=\mathcal{S}_{\sigma_0}(dev)$, where $dev:\widetilde{X}\to \cp$ is 
any developing map for $\sigma$.
\end{definition}

  It is a classical result (see \cite{GU}) that, for any fixed $\sigma_0 \in \mathcal{P}(X)$ 
the map \[ \mathcal{P}(X)\to H^0(X,K^2), \sigma\mapsto \mathcal{S}_{\sigma_0}(\sigma) \] is a 
bijection with the vector space of holomorphic quadratic differentials on $X$, whose dimension is 
$3g-3$ by Riemann-Roch. Of course the zero differential corresponds to the chosen projective 
structure $\sigma_0$; thus the set $\mathcal{P}(X)$ is naturally endowed with the structure of 
complex affine space.

\begin{remark}\label{hownotto2}
From this point of view it is easier to make the observations contained in  Remark  \ref{hownotto1}  clearer. 
 Indeed fix some background $\sigma_0 \in \mathcal{P}(X)$ and pick $\sigma \in \mathcal{P}(X)$; if 
it were possible to consider the quotient $\pi:\sigma \to \sigma/ \psls$ in the category of 
projective structures, then in particular $\pi$ should send the holomorphic quadratic differential 
$q=\mathcal{S}_{\sigma_0}(\sigma)$ on $X$  to a holomorphic quadratic differential on the Riemann 
surface $X/\psls$ corresponding to $\sigma/\psls$. For this to happen we need $\pi$ to branch in a 
controlled way with respect to the zero divisor of  $q$: if $x \in X$ is a point at which the 
stabilizer of the action of $\psls$ has order $m\geq 1$ and at which $q$ has a zero of order $s\geq 
0$, then the induced (meromorphic) quadratic differential on $X/\psls$ will have order 
$\frac{s-2(m-1)}{m}$, so that it is holomorphic when $s\geq 2(m-1)$, but has a genuine pole 
otherwise. For a concrete example consider the case in which $\psls$ is large enough to guarantee 
that $X/\psls$ has genus 0 (e.g. the Fuchsian uniformization for a Hurwitz surface or for a 
hyperelliptic one): in this case the quotient admits no holomorphic quadratic differential 
whatsoever.
\end{remark}

\subsection{Criterion for projectiveness}
  We are now going to consider the action of the biholomorphism group $\aut(X)$ of $X$ on 
the space $\mathcal{P}(X)$ of projective structures on $X$. If $\sigma \in \mathcal{P}(X)$ is 
defined by a developing map $dev:\widetilde{X}\to \cp$ and $F\in \aut(X)$, then $F.\sigma$ is 
defined to be the projective structure with $dev \circ \widetilde{F}^{-1}:\widetilde{X}\to \cp$ as 
a developing map for some lift of $F$ to the universal cover. Since $F$ is holomorphic on $X$, we 
have that $F.\sigma$ is again inside $\mathcal{P}(X)$.\par

Now let us fix a projective structure $\sigma_0 \in \mathcal{P}(X)$. By the above 
discussion we get an identification $\mathcal{P}(X)\cong H^0(X,K^2)$ and we can look at the induced 
action 
on the space of holomorphic differentials, which we denote by $(F,q)\mapsto F.q$. Let us denote by 
$F^*q=q\circ F^{-1}$ the usual action of the automorphism group on the space of differentials by 
pullback. Then a direct computation using the properties of the Schwarzian derivative in  Lemma 
\ref{schwarzianprop} and  Lemma  \ref{schwarzianprop2} shows the following.
\begin{lemma}\label{affineaction}
 Let $F\in \aaut(X)$. If $q=\mathcal{S}_{\sigma_0}(\sigma)$ for some $\sigma \in 
\mathcal{P}(X)$ then $F.q=F^*q+\mathcal{S}_{\sigma_0}(\widetilde{F}^{-1})$.
\end{lemma}
 \begin{proof}
By definition $\mathcal{S}_{\sigma_0}(\sigma)=\mathcal{S}_{\sigma_0}(dev)$ for some developing map 
and $\mathcal{S}_{\sigma_0}(F.\sigma)=\mathcal{S}_{\sigma_0}(dev\circ \widetilde{F}^{-1})$ .
\end{proof}
  In other words $\aut(X)$ acts affinely on $H^0(X,K^2)$, with linear part given by the 
classical 
action by pullback and the translation part accounts for the initial choice of the projective 
structure $\sigma_0$. This action can be used to obtain the following 
criterion.
\begin{proposition}\label{projcriterion}
 Let $F\in \aaut(X)$, $\sigma \in \mathcal{P}(X)$ and $q=\mathcal{S}_{\sigma_0}(\sigma)$. Then 
$F\in \pslss$ if and only if $F.q=q$.
\end{proposition}
\begin{proof}
  Let $dev$ be a developing map for $\sigma$; then $q=\mathcal{S}_{\sigma_0}(dev)$. If 
$F\in 
\psls$ then by  Lemma  \ref{lem_projviadev} we have $dev \circ \widetilde{F}^{-1}=g\circ dev$ for some 
$g \in \pslc$. Therefore $F.q=\mathcal{S}_{\sigma_0}(dev \circ \widetilde{F}^{-1})= 
\mathcal{S}_{\sigma_0}(g\circ dev)= \mathcal{S}_{\sigma_0}(dev)=q$ by  Lemma  \ref{schwarzianprop2}. 
On the other hand if $\mathcal{S}_{\sigma_0}(dev \circ 
\widetilde{F}^{-1})=F.q=q=\mathcal{S}_{\sigma_0}(dev)$ then $dev \circ 
\widetilde{F}^{-1}=g\circ dev$ for some $g \in \pslc$ by  Lemma  \ref{schwarzianprop2}, which implies that 
$F \in \psls$ again by  Lemma  \ref{lem_projviadev}.
\end{proof}

  It is natural to ask if this affine action can be 
reduced to a linear action under a suitable choice of the initial projective structure $\sigma_0$. 
This happens for instance for the Fuchsian uniformization of $X$, as already 
observed in  Lemma  \ref{lem_fuchsunifisrelhurwproj}. 
\begin{corollary}\label{actionislinearforfuchsian}
 Let $\sigma_0$ be the Fuchsian uniformization of $X$, $F\in \aaut(X)$ and $q \in H^0(X,K^2)$. 
Then $F.q=F^*q$.
\end{corollary}
\begin{proof}
Since $\sigma_0$ is the Fuchsian uniformization we have that $F\in \psls$ by  Lemma  
\ref{lem_fuchsunifisrelhurwproj}. Therefore $\mathcal{S}_{\sigma_0}(\widetilde{F}^{-1})=0$ by   Lemma  
\ref{schwarzianprop2}.
\end{proof}
  On the other hand, as shown by the following corollary, there are 
plenty of couples $(\sigma,F)$ where $\sigma$ is a 
projective structure on $X$ and $F$ is a holomorphic automorphism of $X$ which 
is not projective for $\sigma$.

\begin{corollary}
 Let $F\in \aaut(X), F\neq id_X$. If $X$ has genus 2, then also assume $F$ is not the hyperelliptic 
involution. Let $\sigma_0 \in \mathcal{P}(X)$ such that $F\in \pslsspar{\sigma_0}$. Then there 
exists $\sigma \in \mathcal{P}(X)$ such that $F\not \in \pslsspar{\sigma}$.
\end{corollary}

\begin{proof} 
Since $F\in \pslspar{\sigma_0} $, we have that $\mathcal{S}_{\sigma_0}(F)=0$. In particular the 
action of $F$ on $H^0(X,K^2)$ is the linear action by pullback by  Lemma  \ref{affineaction}. The 
linear action of $\aut(X)$ by pullback on $H^0(X,K^2)$ is known to be faithful by 
\cite[V.2]{FK} if and only if $F$ is not the hyperelliptic involution of a genus 2 surface. 
Therefore under our hypothesis there exists $q \in H^0(X,K^2)$ such that 
$F.q=F^*q \neq 
q$. By  Proposition  \ref{projcriterion} we have that $F$ is not projective for the projective structure 
$\sigma$ defined by $\mathcal{S}_{\sigma_0}(\sigma)=q$.
\end{proof}

  For instance recall that by picking the Fuchsian uniformization we can satisfy the 
hypothesis 
of this statement for any $F\in \aut(X)$, by  Lemma  \ref{lem_fuchsunifisrelhurwproj}. We are now ready to 
discuss the following example.
\begin{example}\label{ex_relHurproj}
Let $X$ be a hyperelliptic Riemann surface, and let $J$ be the hyperelliptic involution. By  Proposition 
\ref{projcriterion} and  Corollary  \ref{actionislinearforfuchsian}, the projective structures for which $J$ is 
a projective diffeomorphism are exactly those that correspond to $J$-invariant holomorphic 
quadratic differentials with respect to the Fuchsian uniformization of $X$. Notice that when $g=2$ 
every differential is $J$-invariant, but for $g\geq 3$ the space of $J$-invariant differentials  is 
a proper subspace of $H^0(X,K^2)$, so we get a lot of non-trivial examples of 
non-relatively Hurwitz 
projective structures. On the other hand by choosing $X$ so that $\aut(X)=\{id_X,J\}$ we can obtain 
non-trivial examples of relatively Hurwitz projective structures, namely those corresponding to 
$J$-invariant differentials with respect to the Fuchsian uniformization; 
surfaces of this type exist in any 
genus and were constructed explicitly in \cite{PO}.
  \end{example}

\section{Fuchsian uniformizations and Galois Bely\u{\i} curves}
  In the previous sections we have seen that Fuchsian uniformizations are examples of 
relatively Hurwitz projective structures (\ref{lem_fuchsunifisrelhurwproj}); on the other hand we 
have 
provided lots of examples of relatively Hurwitz projective structures which are not related to 
Fuchsian uniformization (see  Examples  \ref{ex_trivialrelHurproj} and \ref{ex_relHurproj}). Moreover by  Proposition 
\ref{projcriterion} relatively Hurwitz projective structures can be seen as fixed points of an 
affine action of a finite group, so that either there is a unique one, or there is a positive 
dimensional locus of them.
In this section we look for conditions on the underlying Riemann surface under which the Fuchsian 
uniformization is the unique relatively Hurwitz projective structure on it. As it turns out, such a 
condition can be conveniently expressed in terms of  the behaviour of the quotient by the automorphism group.

\begin{definition}
We say that a Riemann surface $X$  is Galois Bely\u{\i}  if $X/\aut(X)$ is 
biholomorphic to $\cp$ and the quotient map $\pi:X\to X/\aut(X)$ is ramified exactly on three 
points.
\end{definition}

The name is motivated by the fact that the quotient map $\pi:X\to X/\aut(X)$ is both a Bely\u{\i} function on $X$ and a (branched) Galois cover. Recall that a Bely\u{\i} function on a Riemann surface is a non-constant meromorphic function ramified over at most three points. Riemann surfaces admitting such a function are called Bely\u{\i} curves, and are known to be exactly those surfaces which can be defined over $\overline{\rat}$ as algebraic curves. 
 
The above condition has several classical equivalents, which we list here for the reader's convenience.  Recall that a $(a,b,c)$-hyperbolic triangle group $\Delta(a,b,c)$ is the orientation-preserving index-two subgroup of the group of isometries $\overline{\Delta}(a,b,c)$ of the hyperbolic plane $\hyp^2$ which is generated by reflections in the sides of a triangle with angles $\frac{\pi}{a},\frac{\pi}{b},\frac{\pi}{c}$. 

\begin{theorem}[Wolfart \cite{WO}]\label{thm_vl} 
Let $X$ be a Riemann surface. Then the following are equivalent:
 \begin{itemize}
  \item[(i)] $X$  is Galois Bely\u{\i} .
  \item[(ii)] $X\cong\hyp^2/\Gamma$, where $\Gamma$ is a cocompact torsion-free 
Fuchsian group whose normalizer in $\pslr$ is a hyperbolic triangle group 
$\Delta$.
  \item[(iii)] $X$ is an isolated local maximum for the function $Y\mapsto \#\aaut(Y)$.
 \end{itemize}
\end{theorem}

\begin{remark}
Such surfaces are also known as quasiplatonic surfaces (terminology introduced by Singerman in \cite{SI}). We refer to \cite[Chapter 5]{JW} for more details and background about their theory.
\end{remark}

The following result allows us to add an item to the previous list in Theorem \ref{thm_vl}.

\begin{theorem}\label{mainthm}
 Let $X$ be a Riemann surface. Then the following are equivalent:
 \begin{itemize}
  \item[(i)] $X$  is Galois Bely\u{\i} .
  \item[(ii)] The Fuchsian uniformization of $X$ is the unique relatively Hurwitz projective 
structure in 
$\mathcal{P}(X)$.
 \end{itemize}
\end{theorem}
\begin{proof}
Let us consider the Schwarzian parametrization of $\mathcal{P}(X)$ with respect to the Fuchsian 
uniformization $\sigma$, so that the action of $\aut(X)$ is linear on $H^0(X,K^2)$ and the zero 
differential corresponds to $\sigma$. By  Proposition  \ref{projcriterion} and  Corollary  \ref{actionislinearforfuchsian} a 
projective structure in $\mathcal{P}(X)$ is relatively Hurwitz if and only if the corresponding 
holomorphic quadratic differential is invariant by this action of $\aut(X)$.
Invariant holomorphic quadratic differentials constitute a linear subspace of $H^0(X,K^2)$, whose 
dimension can be computed according to the following formula (see \cite[V.2.2]{FK})
\[ \dim\left(H^0(X,K^2)^{\aut(X)}\right)=3g_0-3+n \]
where $g_0$ is the genus of $X/\aut(X)$ and $n$ is the number of critical values of the quotient 
map 
$X\to X/\aut(X)$.
Therefore we see that the zero differential is the unique invariant differential exactly 
when $X/\aut(X)$ has genus $0$ and the quotient map is ramified 
exactly over three points (since we always assume the genus of $X$ 
to be at least 2, the case $g_0=1,n=0$ is not allowed).
\end{proof}

\begin{remark}
To get an intuition about the dimension computation in the previous proof, it is worth noticing that, from a Teichm\"uller-theoretic point of view, the Schwarzian parametrization provides an identification between $\mathcal P (X)$ and the cotangent space to Teichm\"uller space at $X$. The space of relatively Hurwitz projective structures corresponds to the cotangent space at $X$ to the locus of fixed points for the action of $\aut(X)$ on Teichm\"uller space, which is a copy of the deformation space for the quotient orbifold $X/ \aut(X)$.
\end{remark}

  By the above theorem this means that complex projective structures which have strictly 
more projective automorphisms than any proper deformations are precisely Fuchsian uniformizations 
of Galois Bely\u{\i} curves. In particular we obtain the following.

\begin{corollary}\label{maincor}
 Let $S$ be a surface and $\sigma$ be a projective structure on it. Then 
$\sigma$ is a Hurwitz projective 
structure if and only if $\sigma$ is the Fuchsian uniformization for a Hurwitz Riemann surface 
structure $X$ on $S$.
\end{corollary}
\begin{proof}
Fuchsian uniformizations of Hurwitz Riemann surfaces are Hurwitz projective by  Lemma 
\ref{lem_fuchshurwishurwproj}. Conversely let $\sigma$ be a Hurwitz projective structure; in 
particular the underlying Riemann surface $X$ is a Hurwitz Riemann surface and $\sigma$ is a 
relatively Hurwitz projective structure on it.
 Since Hurwitz Riemann surfaces are Galois Bely\u{\i} (for instance by Theorem \ref{thm_vl}), we have that Theorem \ref{mainthm} applies  
and implies that $\sigma$ must be the  Fuchsian uniformization of $X$.
\end{proof}

  After establishing such a result, it is natural to ask how often one encounters one of 
the  structures covered by  Theorem  \ref{mainthm} and if something special can be said about its group of 
 projective automorphisms (as done in \cite{SPWS} for the case of translation surfaces). 
 Sporadic and well-studied examples are given by Hurwitz surfaces and Fermat 
 curves $F_n=\{ x^n+y^n=z^n \} \subset \mathbb{CP}^2$; to obtain examples in any genus $g\geq 2$ 
one can consider the hyperelliptic curve defined by the equation $y^2+x^m=1$
 in $\cmp^2$ for all values of $m \in \nat$. On the other hand only  finitely many  Galois Bely\u{\i} curves                                                                                                                          exist for any fixed genus (see \cite{POP} and \cite{SPW}; also see 
\cite[Chapter 5]{JW} for more details and examples). As far as the groups are concerned,  Hurwitz 
groups (i.e.  groups arising as the group of automorphisms of a Hurwitz Riemann  surface) are well-known to be 
quotients of the $(2,3,7)$-triangle group (see \cite{CO} and \cite{CO2}). More generally, as said 
above in  Theorem  \ref{thm_vl}, groups arising as the group of automorphisms of  a Galois Bely\u{\i} curve   is a quotient   of some $(a,b,c)$-hyperbolic triangle group, i.e. they have a 
presentation of the form $\langle x,y \ | \ x^a,y^b,(xy)^c,R\rangle$ for some $a,b,c \in \nat $ 
such that $\dfrac{1}{a}+\dfrac{1}{b}+\dfrac{1}{c}<1$ and for some  extra relation $R$. 

\printbibliography

\end{document}